\documentclass[12pt]{article}
\usepackage[left=2.4cm,top=2.8cm,right=2.4cm,bottom=2.6cm]{geometry}

\usepackage[ansinew]{inputenc}
\usepackage{amscd,amsfonts,amsmath,amssymb,amsthm,enumerate,epsfig,float,graphics,graphicx,pst-all,yhmath}
\newtheorem{lemma}{Lemma}

\newtheorem{theorem}{Theorem}

\newtheorem{remark}{Remark}

\newtheorem{conjecture}{Conjecture}

\newcommand{\fin}{\hfill $\Box$}

\title {An equichordal characterization of the ellipsoid and the sphere}
\author{V. A. Aguilar-Arteaga$^1$, R. I. Ayala-Figueroa$^2$, J. Jer\'onimo-Castro$^3$,\\ and E. Morales-Amaya$^4$\\ 
\small{$^{1,3}$Facultad de Ingenier\'ia}\\
\small{Universidad Aut\'onoma de Quer\'etaro, M\'exico}\\
\small{$^{2}$Instituto Tecnol\'ogico de Mexicali}\\
\small{Tecnol\'ogico Nacional de M\'exico, M\'exico}\\
\small{$^{4}$Facultad de Matem\'aticas-Acapulco,}\\
\small{Universidad Aut\'onoma de Guerrero, M\'exico}\\
\small{$^{1}$aguilar.arteaga@uaq.mx,}\\
\small{$^{2}$rafaelivan@itmexicali.edu.mx,}\\
\small{$^{3}$jesus.jeronimo@uaq.mx,}\\
\small{$^{4}$emoralesamaya@gmail.com}
}
\begin{document}

\maketitle

\begin{abstract}
Let $K$ and $L$ be two convex bodies in $\mathbb R^n$, $n\geq 3$, with $L\subset \text{int}\, K$. In this paper we prove the following result: if every two parallel chords of $K$, supporting $L$ have the same length, then $K$ and $L$ are homothetic and concentric ellipsoids. We also prove a similar theorem when instead of parallel chords we consider concurrent chords. We may also replace, in both theorems, supporting chords of $L$ by supporting sections of constant width. In the last section we also prove similar theorems where we consider projections instead of sections.
\end{abstract}

\section{Introduction}
Let $K$ be a convex body, i.e., a compact and convex set with non-empty interior. We say that a point $x$ in the interior of $K$ is an \emph{equichordal point} if all the chords of $K$ through $x$ have the same length. The famous Equichordal Problem, due to W. Blaschke, W. Rothe, and R. Weitzenb\"ock \cite{Blaschke}, and Fujiwara \cite{Fujiwara}, asks about the existence of a convex planar body with two equichordal points.   M. Rychlik in \cite{Rychlik}, finally gave a complete proof about the non existence of a body with two equichordal points. In \cite{Jero}, the  following extension of
the notion of equichordal point was introduced: Let $K$ and $L$ be two convex bodies in $\mathbb R^n$, $n\geq 2$, with $L\subset \text{int}\, K$; it is said that $L$ is an \emph{equichordal body} for $K$ if every chord of $K$ tangent to $L$ have length equal to a given number $\lambda$. In \cite{Barker}, J. Barker and D. Larman proved that if $K$ is a convex body that admits an equichordal ball then it is also a ball. This result was extended in \cite{Jero} and is proved there that only Euclidean balls possess an equichordal convex body.

Another interesting and classical result is due to W. S\"uss \cite{Suss}: if all $2$-dimensional sections of a convex body $K\subset \mathbb R^3$ through a point $p\in \text{int} K$ are sets of constant width, then $K$ is a ball. This result was latter generalized to $\mathbb R^n$, $n \geq 3$, by L. Montejano in \cite{Montejano}. 

In this paper we give some more characterizations of the Euclidean ball and the ellipsoid: when every set of parallel chords tangent to a convex body in the interior of another convex body have the same length; and when all the sections tangent to the inner convex body, and parallel to a given direction, have constant width.  

\section{A characterization of the ellipsoid}
\begin{figure}[H]
    \centering
    \includegraphics[width=.37\textwidth]{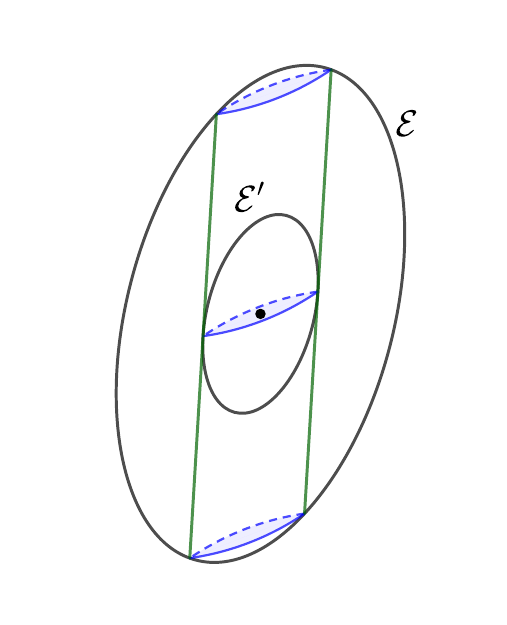}
    \caption{Parallel chords tangent to $\mathcal E'$ have the same length}
    \label{elipsoides}
\end{figure}

Let $\mathcal E\subset \mathbb R^n$, $n\geq 2$, be an ellipsoid and let $\mathcal E'$ be an ellipsoid in its interior, homothetic and concentric with $\mathcal E$. It is not difficult to see that every pair of chords of $\mathcal E$, parallel and tangent to $\mathcal E'$, have the same length (see Fig. \ref{elipsoides}). However, it is very unexpected that this condition characterizes the ellipsoid if the dimension of the space is $n\geq 3$. In the plane, every pair of centrally symmetric and concentric convex bodies share the same property. 

\begin{theorem}\label{parallel}
Let $K,L \subset \mathbb{R}^{n}$, $n\geq 3$, be convex bodies with $L \subset \emph{int} K$. Suppose that for every $u \in \mathbb{S}^{n-1}$ all the chords of $K$ supporting $L$ and parallel to $u$, have the same length $\lambda (u)$. Then $K$ and $L$ are homothetic and concentric ellipsoids.
\end{theorem}

\emph{Proof.} Let $u \in \mathbb{S}^{n-1}$ be any unit vector. Let $[a,b]$ be any chord of $K$ parallel to $u$ touching $L$ at $x$. Consider a $2$-dimensional plane $H$ supporting $L$ at $x$ and such that $[a,b] \subset H$. We will show that $[a,b]$ is an affine diameter of $H \cap K$. That is, there exist two parallel supporting $(n-2)$-dimensional planes of $H \cap K$ at points $a$ and $b$, or equivalently, $[a,b]$ is a longest chord of $H \cap K$ parallel to $u$. Suppose this is not the case and let $[a',b']$ be the affine diameter of $H \cap K$ parallel to $[a,b]$. We have that $|a'b'|>|ab|$. Let $\Pi$ be a $2$-dimensional plane containing $[a'b']$ such that $\Pi$ intersect the interior of $L$. Denote by $[c,d]$ and $[c',d']$ the two chords of $K$, parallel to $[a',b']$, supporting $\Pi \cap L$. By hypothesis $|cd|=|c'd'|<|a'b'|$, and by the convexity of $\Pi \cap K$, this is not possible. Hence $|ab|=|a'b'|$ and then $[a,b]$ is an affine diameter of $H \cap K$. In the same way we can prove that any other chord of $H \cap K$ through $x$ is an affine diameter of $H \cap K$. By a well known result of Hammer \cite{Hammer} we have that $H \cap K$ has center of symmetry at the point $x$. This is true for every $2$-dimensional plane supporting $L$ at $x$, it follows that the hypersection supporting $L$ at $x$ has center of symmetry at $x$. Now, by a theorem due to Olovjanishnikov \cite{Olovja} we have that $K$ and $L$ are homothetic and concentric ellipsoids. \fin

We believe the following conjectures are true.

\begin{conjecture}
Let $K,L \subset \mathbb R^2$ be convex bodies with $L\subset \emph{int} K$ and $L$ centrally symmetric. Suppose every pair of parallel chords of $K$ supporting $L$ have the same length. Then $K$ is also centrally symmetric.
\end{conjecture}

\begin{conjecture}
Let $K,L \subset \mathbb R^3$ be convex bodies with $L\subset \emph{int} K$ a strictly convex body. Suppose every section of $K$ supporting $L$ has the contact point as an equichordal point. Then $K$ and $L$ are concentric balls.
\end{conjecture}


\section{Characterization of the sphere by concurrent chords}
The following lemma is interesting by itself and will be useful for the subsequent results.

\begin{lemma}\label{la_elipse}
Let $\mathcal E_1$ and $\mathcal E_2$ be two concentric and homothetic ellipses with $\mathcal E_1\subset \emph{int} \mathcal E_2$. Let $[a,a']$ be the chord of $\mathcal E_2$, tangent to $\mathcal E_1$, and orthogonal to the mayor axis of $\mathcal E_2$. Then, any other chord $[b,b']$ tangent to $\mathcal E_1$, not orthogonal to the mayor axis of $\mathcal E_2$, has length bigger than the length of $[a,a']$.
\end{lemma}

\begin{figure}[H]
    \centering
    \includegraphics[width=.53\textwidth]{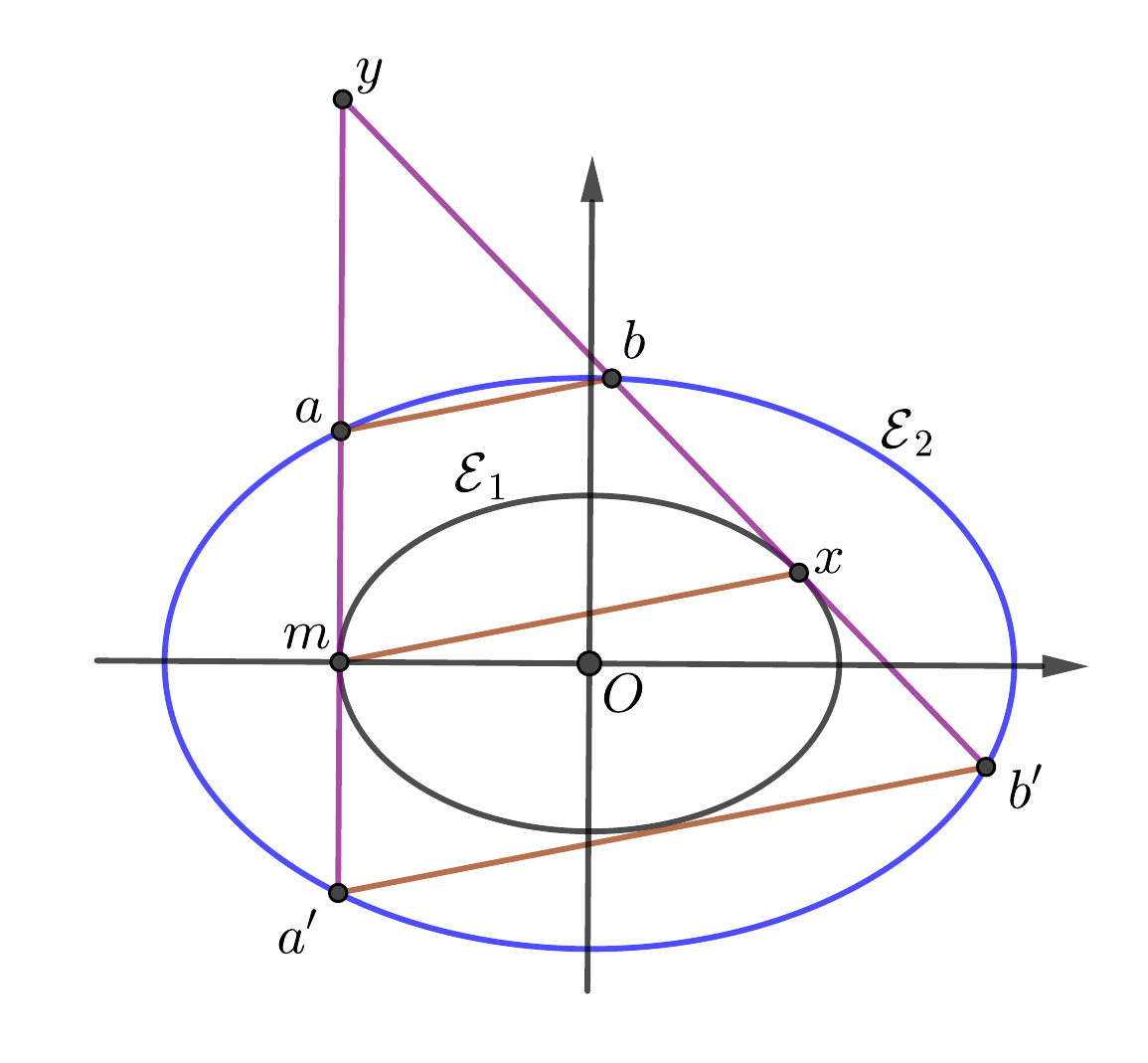}
    \caption{$[a,a']$ is the tangent chord with minimum length}
    \label{elipses}
\end{figure}

\emph{Proof.} Let $y$ be the point of intersection between the lines $aa'$ and $bb'$, and let $m$ and $x$ be the points where the chords $[a,a']$ and $[b,b']$ touch $\mathcal E_1$, respectively. We know that $[a,b] \parallel [m,x] \parallel [a',b']$ (see Fig. \ref{elipses}), this can be seen easily if we apply an affine transformation which send $\mathcal E_1$ and $\mathcal E_2$ to concentric circles. To see that $|aa'|<|bb'|$ it is enough, by Thales's theorem, to prove that $|ym|<|yx|$.

\begin{figure}[H]
    \centering
    \includegraphics[width=.5\textwidth]{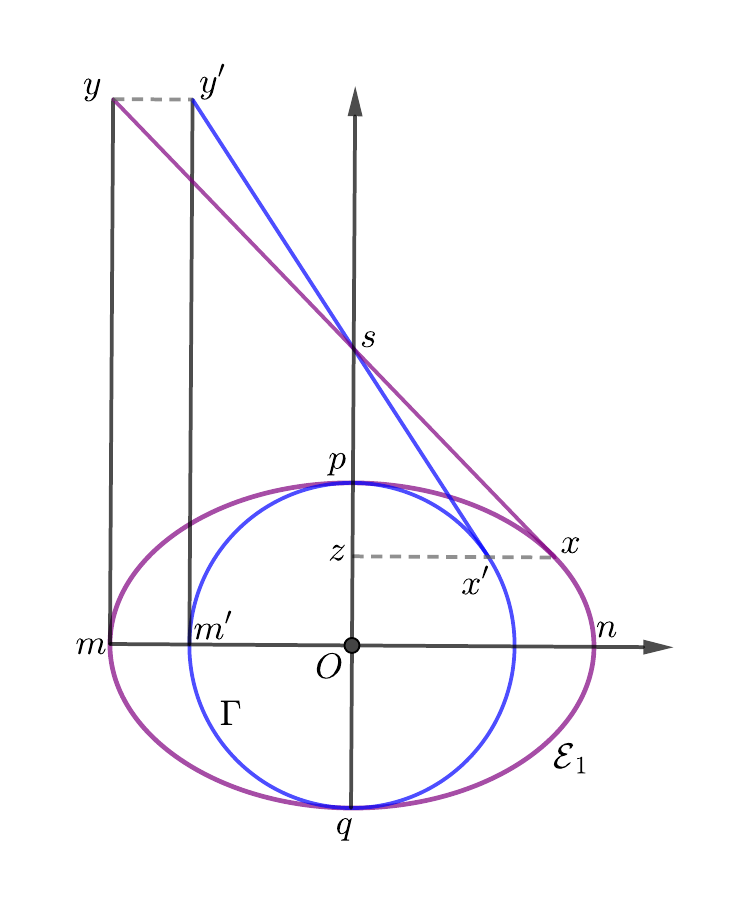}
    \caption{The length of $[y,x]$ is bigger than the length of $[y',x']$}
    \label{elipse}
\end{figure}

Consider the circle $\Gamma$ with diameter $[p,q]$, the minor axis of $\mathcal E_1$ (see Fig. \ref{elipse}). Let $z$ be the projection of $x$ into the minor axis of $\mathcal E_1$ and let $x'$ be the point where the segment $[z,x]$ intersects the circle $\Gamma$. By a well known property of the ellipse, we have that $\frac{|zx|}{|zx'|}=\lambda >1$, for a fixed number $\lambda$. The linear transformation that sends $\mathcal E_1$ to $\Gamma$, also sends the tangent segments $[y,m]$ and $[y,x]$ to the segments $[y',m']$ and $[y',x']$, tangent to $\Gamma$. Clearly, we have that $|ym|=|y'm'|$ and $|y'x'|<|yx|$ and since $|y'm'|=|y'x'|$ we have that $|ym|<|yx|$. \fin

\begin{theorem}\label{concurrent}
Let $K,L, M \subset \mathbb{R}^{n}$, $n\geq 3$, be convex bodies with $L \subset \emph{int} K\subset \emph{int} M$. Suppose that for every $x \in \partial M$ all the lines supporting $L$ and passing through $x$, intersects $K$ in chords of the same length $\lambda (x)$. Then $K$ and $L$ are concentric balls.
\end{theorem}

\emph{Proof.} First we prove the theorem in dimension $n=3$. Let $x \in \partial M$ be any point and $C(x,L)$ be the support cone of $L$ with apex at $x$. Define as $\gamma(x)=\{C(x,L) \cap \partial M\}\setminus \{x\}$, which is clearly a simple and closed curve in the boundary of $M$. Notice that $\gamma(x)$ divides the boundary of $M$ into two open regions, the region $R^-$ that contains $x$ and the region $R^+$ that does not contain $x$. Now, consider a chord $[x,y]$ so that $y \in \gamma(x)$ and any chord $[x',y']$ of $M$ parallel to $[x,y]$, supporting $L$, and such that $x'\in R^-$ and $y'\in R^+$. Since the cone $C(x',L)$ intersects $R^-$ we have that $\gamma(x) \cap \gamma(x') \neq \emptyset$. Choose a point $z \in \gamma(x) \cap \gamma(x')$. Since $zx$ and $zx'$ are support lines of $L$ that pass through $z \in \partial M$, then $\lambda(z)=|[z,x']\cap K|=|[z,x]\cap K|=\lambda(x)=\lambda (x')$. We have proved that the intersection of $[x,y]$ with $K$ and the intersection of $[x',y']$ with $K$ have the same length. That is, parallel chords of $K$ supporting $L$ have the same length. From Theorem \ref{parallel} we can say that $K$ and $L$ are homothetic and concentric ellipsoids. But from Lemma \ref{la_elipse} we can see that $K$ and $L$ must be, indeed, concentric balls. 

For dimension $n>3$ we just use the fact that if every $3$-dimensional section of a convex body is a $3$-dimensional ball, then it is an $n$-dimensional ball. \fin

\begin{remark}
We may replace the set $M$ in Theorem \ref{concurrent} by a pair of parallel hyperplanes $H_1$ and $H_2$ such that $K$ and $L$ are contained in the open region between those planes. The proof follows the same ideas.
\end{remark}


\section{Characterization of the sphere by sections of constant width}

\begin{theorem}\label{parallel_sections}
Let $K,L \subset \mathbb{R}^{3}$ be convex bodies with $L \subset \emph{int} K$. Suppose that for every $u \in \mathbb{S}^{2}$ all sections of $K$ supporting $L$ and parallel to $u$, have the same constant width $h(u)$. Then $K$ and $L$ are concentric balls.
\end{theorem}

\emph{Proof.} We will prove first that $K$ is strictly convex. Suppose this is not the case and let $[p,q]$ be a segment in the boundary if $K$. Let $S_1$ be a supporting $2$-dimensional plane of $L$ containing $[p,q]$. Since $S_1\cap K$ is of constant width, it must be strictly convex, then $[p,q]$ must be a point. Now, let $u \in \mathbb{S}^{2}$ be any unit vector. Let $[a,b]$ be a chord of $K$, parallel to $u$, with length $h(u)$. We will show that $[a,b]$ is a supporting chord of $L$. Let $[x,y]$ be another chord of $K$ with length $h(u)$ such that the parallelogram $abyx$ intersects the interior of $L$. Let $\Pi$ be the plane containing the chords $[a,b]$ and $[x,y]$. Consider a chord $[c,d]\subset \Pi\cap K$, parallel to $u$, tangent to $L$, and different from $[a,b]$ and $[x,y]$. There are two possible positions for $[c,d]$. If $[c,d]$ lies between $[a,b]$ and $[x,y]$, then, since $\Pi\cap K$ is strictly convex, we have $|cd|>h(u)$, which contradicts the hypothesis that every section of $K$ supporting $L$ and parallel to $u$ has constant width $h(u)$. If $[c,d]$ lies outside the strip determined by $[a,b]$ and $[x,y]$, then $|cd|<h(u)$, which again contradicts the same hypothesis. Therefore, no such chord $[c,d]$ can exist, and hence $[a,b]$ must be tangent to $L$. Now, let $\Gamma$ be any plane tangent to $L$ and let $[a,b]$ be the chord of $K$ in $\Gamma$ parallel to $u$ which is tangent to $L$. If $|ab|< h(u)$ then there is a chord of $K$ in $\Gamma\cap K$, parallel to $u$ with length $h(u)$ which is not tangent to $L$. This is a contradiction to the fact proved before. With this, we have proved that all the chords of $K$, parallel to $u$ and supporting $L$, have the same length $h(u)$. This is true for every $u\in\mathbb S^2$; hence, we can apply Theorem  \ref{parallel} and obtain that K and L are homothetic and concentric ellipsoids. However, the only ellipsoids with sections of constant width are Euclidean balls. \fin

In a very similar way we can prove the following.

\begin{theorem}\label{concurrent_sections}
Let $K,L, M \subset \mathbb{R}^{3}$ be convex bodies with $L \subset \emph{int} K\subset \emph{int} M$. Suppose for every $x \in \partial M$ all sections of $K$ supporting $L$ and passing through $x$, are of constant width $h(x)$. Then $K$ and $L$ are concentric balls.
\end{theorem}

\emph{Proof.} Following the same ideas as in \ref{concurrent}, we show that any two parallel sections are of the same constant width. Hence using Theorem \ref{parallel_sections} and Lemma \ref{la_elipse} the result follows directly.
\fin


\section{A  particular case of S\"uss theorem}
We start this section with a proof for a particular case of S\"uss theorem, i.e., we will prove that a convex body with concurrent sections of constant width 1 is a ball. Before giving the proof we introduce the following notation: for
every $\upsilon\in \mathbb S^2$ denote by $S(\upsilon)$ the circle
$\mathbb S^2\cap \upsilon ^{\bot},$ by $H(\upsilon)$ and
$H(-\upsilon)$ the two supporting planes of $K$ orthogonal to
$\upsilon,$ and by $\text{Sb}(K,\upsilon)$ the shadow boundary of $K$ in
direction $\upsilon$, i.e., $$\text{Sb} (K,\upsilon)\equiv \{x\in\partial
K:\text{there is a line parallel to }\upsilon \text{ and touching
} K \text{ at } x. \}$$ It is a well known result of W. 
Blaschke \cite{Blaschke2} that a convex body $K$ is an ellipsoid if for
every $\upsilon\in\mathbb S^2$ the shadow boundary
$Sb(K,\upsilon)$ is a planar curve. Moreover, if we know that
$\text{Sb}(K,\upsilon)$ is orthogonal to $\upsilon$, for every
$\upsilon\in\mathbb S^2$, then the conclusion is even stronger:
$K$ is a Euclidean ball. However, if we restrict the directions for the shadow boundaries we have the following interesting conclusion.

\begin{lemma}{\label{lemma2}} Let $K\subset\mathbb R^3$ be a convex body and $\upsilon\in\mathbb
S^2$ such that for every $w\in S(\upsilon)$ we have that $\emph{Sb}(K,w)$
is a closed planar curve orthogonal to $w$. Then $K$ has an axis
of revolution parallel to $\upsilon$.
\end{lemma}

\emph{Proof.} Let $H(\upsilon)$ and
$H(-\upsilon)$ be the planes supporting $K$, orthogonal to
$\upsilon,$ and consider the points $a\equiv H(\upsilon)\cap K$
and $b\equiv H(-\upsilon)\cap K.$ Clearly, $a$ and $b$ belong to
$\text{Sb}(K,w)$ for every $w\in S(\upsilon)$, and since $w\bot \
\text{Sb}(K,w)$ then we get that $[a,b]$ is parallel to $\upsilon.$ Now, let
$\Omega$ be a plane parallel to $H(\upsilon)$ which intersects the
interior of $K$. Consider an arbitrary vector $w\in S(\upsilon),$
let $\Pi$ be the plane orthogonal to $w$ through $[a,b]$, and let
$\{x,y\}\equiv\Omega\cap\partial K\cap\Pi$. We have that
$\text{Sb}(K,w)=\Pi\cap\partial K$, then the tangent vectors to
$\Omega\cap\partial K$ through the points $x$ and $y$ are
parallel to $w$. We have that $[x,y]$ is orthogonal to $w$ and intersects the segment
$[a,b]$ in a point $m$. It is known that a planar curve is a circle
if all its normal lines are concurrent, it follows that
$\Omega\cap\partial K$ is a circle with center at $m$.
Since the above is true for every plane $\Omega $ orthogonal to
$\upsilon$ and intersecting $K$, we conclude that the line $ab$ is
an axis of revolution for $K$. \fin

\begin{figure}[H]
    \centering
    \includegraphics[width=.68\textwidth]{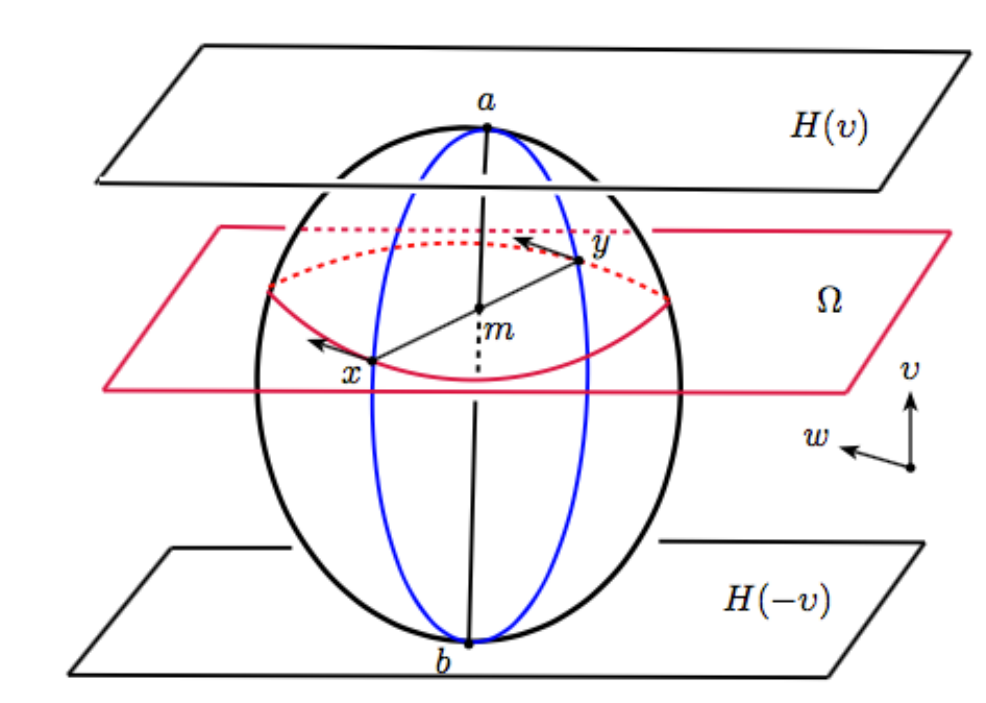}
    \caption{The shadow boundaries orthogonal to $\upsilon$ are planar curves}
    \label{planas}
\end{figure}

\begin{theorem}\label{Suss} Let $K\subset \mathbb R^3$ be a convex body and let $p$ be a point in its interior. If every $2$-dimensional section of $K$ through $p$ has constant width 1, then $K$ is a ball with center at $p$ and diameter 1.
\end{theorem}

\emph{Proof.} Let $[a,b]$ be any affine diameter of $K$. The section containing $p$, $a$, and $b$, has constant width and has diameter equal to $|ab|=1.$ It follows that all the affine diameters of $K$ are indeed diameters and hence $K$ is a body of constant width 1.   Let $H$ be any $2$-dimensional plane containing $p$ with orthogonal unit vector $u$. The orthogonal projection of $K$ onto $H$, $\pi_u(K)$, is a $2$-dimensional body of constant width 1 which contains $H\cap K$, then we have that $\pi_u(K)=H\cap K$. It follows that $\text{Sb} (K,u)=\partial (H\cap K)$. We have proved that all the shadow boundaries of $K$ are planar curves, by Blaschke's theorem we have  that $K$ must be an ellipsoid. However, since $K$ is a body of constant width it must be a ball. \fin

\begin{remark}
This theorem can also be proved using Lemma \ref{lemma2}: it is not difficult to prove that there are two diameters of $K$ passing through $p$. By Lemma \ref{lemma2} we have that each one of these diameters is an axis of revolution of $K$ and is well known that the only convex body with two axes of revolution is the Euclidean ball.
\end{remark}

\section{Equichordal projections}

We give a result where instead of sections tangent to the inner body we consider orthogonal projections.

\begin{theorem}\label{proyeccion_equicordal}
Let $K,L \subset \mathbb{R}^{3}$ be convex bodies with $L \subset \emph{int} K$ and $L$ strictly convex. Suppose that for every $u \in \mathbb{S}^2$ all the chords of $\pi_u(K)$ that are tangent to $\pi_u(L)$ have length equal to 1. Then $K$ and $L$ are concentric balls.
\end{theorem}

\begin{figure}[H]
    \centering
    \includegraphics[width=1.0\textwidth]{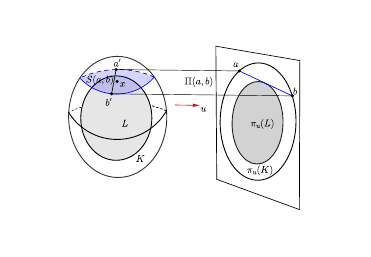}
    \caption{The section $S(a,b)$ is a set of constant width $1$}
    \label{proyeccion}
\end{figure}

\emph{Proof.} Let $u\in\mathbb S^2$ be any direction and let $[a,b]$ be any chord of $\pi_u(K)$ supporting $\pi_u(L)$. Let $\Pi(a,b)$ be the plane parallel to $u$ and containing $[a,b]$. Clearly, $\Pi(a,b)$ is a supporting plane of $L$. Let $a', b'\in \partial K$ be the points such that $\pi_u(a')=a$ \ and $\pi_u(b')=b.$ We have that the lines $a'a$ and $b'b$ are supporting lines of the section $S(a,b)=\Pi(a,b)\cap K$ (see Fig. \ref{proyeccion}). If $[a',b']$ is not parallel to $[a,b]$ then we can find a direction $w\in\mathbb S^2$ where a chord of $\pi_w(K)$ tangent to $\pi_w(L)$ has length bigger than $1$. To see this, consider $w$, orthogonal to $[a',b']$ and parallel to $\Pi(a,b)$. We have that the segment $\Pi(a,b)\cap \pi_w(K)$ has length bigger than or equal to $|a'b'|>|ab|=1$, which is a contradiction to the hypothesis of the theorem. Hence we have that $[a',b']$ is parallel to $[a,b]$, i.e., the width of $S(a,b)$ in the direction of $[a,b]$ is equal to $1$. In the same way we prove that the width of $S(a,b)$ in every direction parallel to $\Pi(a,b)$ is equal to $1$. In other words, $S(a,b)$ is a set of constant width $1$. We have proved that all the sections of $K$ supporting $L$ are sets of constant width $1$, we apply Theorem \ref{parallel_sections} and conclude that $K$ and $L$ are concentric balls. \fin

\begin{remark}
Following similar ideas, we can prove the case in Theorem \ref{proyeccion_equicordal} when $L$ is a point. However, we suspect the following is also true.
\end{remark}

\begin{conjecture}
Let $K,L \subset \mathbb{R}^{3}$ be convex bodies with $L \subset \emph{int} K$ and $L$ strictly convex. Suppose that for every $w \in \mathbb{S}^2$ the chords of $\pi_u(K)$ that are tangent to $\pi_u(L)$ and are parallel to $w$, with $u\bot w$, have length equal to $\lambda (w)$. Then $K$ and $L$ are homothetic and concentric ellipsoids.
\end{conjecture}

Finally, we have the following conjecture for which in the last theorem we provide a partial result towards the proof of the full conjecture.

\begin{conjecture}
Let $K \subset \mathbb{R}^{3}$ be a convex body and let $p$ be a point in its interior. Suppose that for every $u \in \mathbb{S}^2$ the point $\pi_u(p)$ is an equichordal point of $\pi_u(K)$. Then $K$ is ball.
\end{conjecture}

\begin{theorem}\label{punto_equicordal}
Let $K\subset \mathbb{R}^3$ be a strictly convex body and let $p$ be an equichordal point of $K$. Suppose that for all $u\in \mathbb{S}^2$, $\pi_u(p)$ is an equichordal point of $\pi_u(K)$. Then $K$ is a body of revolution.
\end{theorem}

\emph{Proof.} Let $[a,b]$ be an affine diameter of $K$ through $p$. We will prove now that $[a,b]$ is a binormal of $K$, i.e., there exist support planes of $K$ through $a$ and $b$, respectively, which are orthogonal to the chord $[a,b]$. Let $\Pi_a$ and $\Pi_b$ be parallel support planes of $K$ through $a$ and $b$, respectively, and let $v\in\mathbb S^2$ be the unit vector orthogonal to them. If $v$ is parallel to $[a,b]$, then $[a,b]$ is indeed a binormal of $K$ and we are done. Let us assume that $v$ is not parallel to $[a,b]$. Let $w\in\mathbb S^2$ be parallel to $[a,b]$ and let $u$ be a unit vector orthogonal to $v$ and parallel to the subspace generated by $v$ and $w$. Let $a'=\pi_u(a)$ and $b'=\pi_u(b)$ be the orthogonal projections of $a$ and $b$ onto the plane $u^{\bot}$. Let $[x,y]$ be the chord of $K$ through $p$ which is parallel to $u^{\bot}$ and consider the segment $[x',y']\subset \pi_u(K)$, with $x'=\pi_u(x)$ and $y'=\pi_u(y)$. Since $[x',y']$ passes through $p'=\pi_u(p),$ which is an equichordal point of $\pi_u(K)$, we have that $$|a'b'|\geq |x'y'|=|xy|=|ab|,$$ however, this is only possible if $|a'b'|=|ab|$, which implies that $[a,b]$ is parallel to $[a',b']$. We have proved that $[a,b]$ is orthogonal to $\Pi_a$ and $\Pi_b$, i.e., $[a,b]$ is a binormal of $K$.

Now, let $u\in \mathbb{S}^2$ be any vector orthogonal to $[a,b]$ and denote by $K_u$ the $2$-dimensional section of $K$ orthogonal to $u$ and passing through $p$. We have that all chords of $\pi_u(K)$ through $\pi_u(p)$ have length equal to $|\pi_u(a)\pi_u(b)|=|ab|$, and since $\pi_u(K_u)\subset \pi_u(K)$ we have that $\pi_u(K_u)=\pi_u(K)$. Since $K$ is a strictly convex body, we obtain that the shadow boundary of $K$ in direction $u$, $\text{Sb} (K,u)$, coincides with $\partial K_u$, in other words, $\text{Sb}(K,u)$ is a planar closed curve orthogonal to $u$. Now we apply Lemma \ref{lemma2} and conclude that the line $ab$ is an axis of revolution for $K$. \fin


\section{Statements and Declarations}

\textbf{Funding}

The authors declare that no funds, grants, or other support were received during the preparation of this manuscript.

\textbf{Competing Interests}

The authors have no relevant financial or non-financial interests to disclose.

\textbf{Author Contributions}

All authors contributed to the study conception and design. Material preparation were performed by Jes\'us Jer\'onimo Castro, Rafael Ayala Figueroa, Efr\'en Morales Amaya, and V\'ictor Aguilar Arteaga. The first draft of the manuscript was written by Jes\'us Jer\'onimo and Rafael Ayala and all authors commented on previous versions of the manuscript. All authors read and approved the final manuscript.

\textbf{Data availability}

Data sharing not applicable to this article as no datasets were generated or analysed during the current study.


\end{document}